\documentclass[12pt,a4paper,oneside,reqno]{article}

\usepackage[utf8]{inputenc}
\usepackage[english]{babel}
\usepackage[left=1in,right=1in,top=1in,bottom=1in]{geometry}
\usepackage{amsfonts,amssymb,amsmath,amsthm}
\usepackage{xcolor}
\usepackage[normalem]{ulem}
\usepackage{cite}
\usepackage{hyperref}
\hypersetup{colorlinks=true,linkcolor=blue,citecolor=blue,urlcolor=blue}
\usepackage{fancyhdr}
\newtheorem{theorem}{Theorem}
\newtheorem{lemma}[theorem]{Lemma}

\newcommand{\doi}[1]{\href{https://doi.org/#1}{doi:#1}}
\newcommand{\arxiv}[1]{\href{https://arxiv.org/abs/#1}{arXiv:#1}}

\begin{document}

\noindent\fbox{\parbox{\dimexpr\textwidth-2\fboxsep-2\fboxrule\relax}{\small
\textbf{Corrected version verified in Lean.} This is a corrected version of the paper
I.~Kaygorodov, \emph{On $\delta$-derivations of $n$-ary algebras}, Izv. Math. 76 (2012),
no.~6, 1150--1162, \doi{10.1070/IM2012v076n06ABEH002618}.}}

\bigskip

\begin{center}
{\Large On $\delta$-derivations of $n$-ary algebras}

\bigskip

{\bf Ivan Kaygorodov}

\end{center}

\medskip

\noindent{\small\textbf{Abstract.} We give a description of $\delta$-derivations of
$(n+1)$-dimensional $n$-ary Filippov algebras and, as a consequence, of simple
finite-dimensional Filippov algebras over an algebraically closed field of characteristic
zero. We also give new examples of non-trivial $\delta$-derivations of Filippov algebras and
show that there are no non-trivial $\delta$-derivations of the simple ternary Malcev 
algebra $M_8$.}

\medskip

\noindent{\small\textbf{Keywords:} $\delta$-derivation, Filippov algebra, ternary Malcev 
algebra.}

\bigskip

\noindent {\bf MSC2020:} 17A36 (primary); 17A40, 17A42 (secondary).

\bigskip

\setcounter{tocdepth}{1}
\tableofcontents

\section{Introduction}

The notion of anti-derivation of an algebra, which is a special case of a
$\delta$-derivation, that is, a $(-1)$-derivation, was considered in~\cite{b1,b2}.
Subsequently the definition of a $\delta$-derivation of an algebra appeared in~\cite{b3}.
Recall that for a fixed $\delta$ in the ground field $F$, a $\delta$-derivation of an
algebra $A$ is a linear map $\varphi$ satisfying the condition
\begin{equation}\label{e1}
\varphi(xy)=\delta\bigl(\varphi(x)y+x\varphi(y)\bigr)
\end{equation}
for arbitrary elements $x,y\in A$.

A description of $\frac12$-derivations of an arbitrary prime Lie algebra $A$
$\bigl(\frac16\in F\bigr)$ with a non-degenerate symmetric invariant bilinear form was given
in~\cite{b3}: it was proved that a linear map $\varphi\colon A\to A$ is a
$\frac12$-derivation if and only if $\varphi\in\Gamma(A)$, where $\Gamma(A)$ is the centroid
of $A$. Hence, if $A$ is a central simple Lie algebra over a field of characteristic
$p\neq2,3$ with a non-degenerate symmetric invariant bilinear form, then any
$\frac12$-derivation $\varphi$ has the form $\varphi(x)=\alpha x$ for some $\alpha\in F$.
Filippov~\cite{b4} proved that a prime Lie algebra cannot have a non-zero
$\delta$-derivation if $\delta\neq-1,0,\frac12,1$. It was shown in~\cite{b4} that a prime Lie
algebra $A$ $\bigl(\frac16\in\Phi\bigr)$ with a non-zero anti-derivation is a 3-dimensional
central simple algebra over the field of fractions of the centre $Z_R(A)$ of its algebra of
right multiplications $R(A)$. In the same paper, an example was constructed of a
non-trivial $\frac12$-derivation of the Witt algebra $W_1$, that is, a $\frac12$-derivation
that is not an element of the centroid of $W_1$. In~\cite{b5}, $\delta$-derivations were
described of prime alternative and non-Lie Malcev  algebras with certain restrictions on
the ring of operators $\Phi$. It turned out that algebras in these classes have no
non-trivial $\delta$-derivations.

A description of $\delta$-derivations of simple finite-dimensional Jordan superalgebras
over an algebraically closed field of characteristic zero was given in~\cite{b6}, and
$\delta$-derivations of classical Lie superalgebras were described in~\cite{b7}. The
paper~\cite{b8} was devoted to a description of $\delta$-derivations of semisimple
finite-dimensional Jordan superalgebras over an arbitrary field of characteristic different
from~$2$ and of $\delta$-(super)derivations of simple finite-dimensional Lie and Jordan
superalgebras over an algebraically closed field of characteristic zero. It was shown that
the algebras and superalgebras in~\cite{b6,b7,b8} have no non-trivial
$\delta$-(super)derivations. Subsequently the results in~\cite{b7} were generalized by
Zusmanovich~\cite{b9}, who gave a description of $\delta$-(super)derivations of prime Lie
superalgebras, namely, he proved that a prime Lie superalgebra has no non-trivial
$\delta$-(super)derivations for $\delta\neq-1,0,\frac12,1$. He also showed that if a Lie
superalgebra $A$ is perfect ($[A,A]=A$) and has zero centre and a non-degenerate
supersymmetric invariant bilinear form, then the space of $\frac12$-(super)derivations
coincides with the (super)centroid of $A$. Furthermore, in the case of positive
characteristic he gave an affirmative answer to Filippov's question in~\cite{b4} on the
existence of zero-divisors in the ring of $\frac12$-derivations of a prime Lie algebra. The
present author and Zhelyabin~\cite{b10} considered $\delta$-(super)derivations of simple
unital superalgebras of Jordan brackets. They showed that there are no non-trivial
$\delta$-(super)derivations of simple superalgebras of Jordan brackets that are not
superalgebras of vector type and gave a description of the $\delta$-(super)derivations of
simple finite-dimensional unital Jordan superalgebras over an algebraically closed field of
characteristic different from~$2$. As a consequence, a connection was discovered between the
existence of non-trivial $\delta$-derivations of simple unital superalgebras of Jordan
brackets and the speciality of the superalgebra. The cycle of papers describing the
$\delta$-(super)derivations of simple finite-dimensional Jordan superalgebras includes the
paper~\cite{b11}, which contains a complete description of the $\delta$-(super)derivations
of semisimple finite-dimensional Jordan superalgebras over an algebraically closed field of
characteristic different from~$2$. In particular, examples of non-trivial
$\frac12$-(super)derivations were constructed for simple non-unital finite-dimensional
Jordan superalgebras.

Here we consider $\delta$-derivations of $n$-ary Filippov algebras, as well as of the
ternary Malcev  algebra $M_8$. We point out that one of the sources of origin of Filippov
algebras was the Nambu mechanics proposed in~\cite{b12}. Broad connections are currently
known between problems in the theory of $n$-ary algebras and mathematical physics, in
particular, between Nambu mechanics and Chern--Simons theory (for more details,
see~\cite{b13,b14,b15,b16,b17,b18,b19}).

\section{New examples of non-trivial $\delta$-derivations of Filippov algebras}

A Filippov algebra is an algebra $L$ with one anticommutative $n$-ary operation
$[x_1,\dots,x_n]$ satisfying the identity
\[
\bigl[[x_1,\dots,x_n],y_2,\dots,y_n\bigr]=\sum_{i=1}^n\bigl[x_1,\dots,[x_i,y_2,\dots,y_n],\dots,x_n\bigr].
\]

A derivation of an $n$-ary algebra $L$ is a linear map $D$ satisfying the condition
\begin{equation}\label{e2}
D[x_1,\dots,x_n]=\sum_{i=1}^n[x_1,\dots,D(x_i),\dots,x_n].
\end{equation}
By analogy with $\delta$-derivations of binary algebras, as studied in the
papers~\cite{b3,b4,b5,b6,b7,b8,b9,b10,b11}, for a fixed element $\delta$ of the ground field
we can define a $\delta$-derivation of an $n$-ary algebra as a linear map $\varphi$
satisfying the condition
\begin{equation}\label{e3}
\varphi[x_1,\dots,x_n]=\delta\sum_{i=1}^n[x_1,\dots,\varphi(x_i),\dots,x_n].
\end{equation}

Let
\[
\Gamma(L)=\{\psi\in\mathrm{End}(L)\mid
\psi([x_1,\dots,x_n])=[x_1,\dots,\psi(x_i),\dots,x_n]\}
\]
be the centroid of an algebra $L$. Clearly, in the case of an $n$-ary algebra $L$, every
element of $\Gamma(L)$ is a $\frac1n$-derivation. A non-zero $\delta$-derivation $\varphi$ is
considered to be non-trivial if $\delta\neq0,1$ and $\varphi\notin\Gamma(L)$. The cases
$\delta=0,1$ are trivial and therefore are not considered in what follows.

Filippov~\cite{b20} gave a description of $n$-ary algebras, which were subsequently called
Filippov algebras, of dimension $n+1$. Let $A$ be an $n$-ary Filippov algebra,
$\dim(A)=n+1$, and let $\{e_1,\dots,e_{n+1}\}$ be a basis of $A$. We write $\hat x_i$ to
denote the absence of the component $x_i$. For example,
$[x_1,x_2,\hat x_3,x_4]=[x_1,x_2,x_4]$.

In accordance with Filippov's classification, the $n$-ary Filippov algebras of dimension
$n+1$ over an algebraically closed field $F$ of characteristic different from~$2$ are
exhausted by the following $n+4$ families of algebras:
\begin{itemize}
\item[$(A_1)$] $[e_1,\dots,\hat e_i,\dots,e_{n+1}]=0$ for $1\leqslant i\leqslant n+1$;
\item[$(B_1)$] $[e_1,\dots,\hat e_i,\dots,e_{n+1}]=0$ for $1<i\leqslant n+1$,
$[e_2,\dots,e_{n+1}]=e_1$;
\item[$(B_2)$] $[e_1,\dots,\hat e_i,\dots,e_{n+1}]=0$ for $1\leqslant i\leqslant n$,
$[e_1,\dots,e_n]=e_1$;
\item[$(C_1)$] $[e_1,\dots,\hat e_i,\dots,e_{n+1}]=0$ for $1\leqslant i\leqslant n-1$,
$[e_1,\dots,\hat e_n,e_{n+1}]=e_n$, $[e_1,\dots,e_n]=\alpha e_{n+1}$, $\alpha\neq0$;
\item[$(C_2)$] $[e_1,\dots,\hat e_i,\dots,e_{n+1}]=0$ for $1\leqslant i\leqslant n-1$,
$[e_1,\dots,\hat e_n,e_{n+1}]=e_n+\beta e_{n+1}$, $[e_1,\dots,e_n]=e_{n+1}$, $\beta\neq0$;
\item[$(D_r)$] $[e_1,\dots,\hat e_i,\dots,e_{n+1}]=e_i$ for $3\leqslant r\leqslant n+1$ and
$1\leqslant i\leqslant r$, $[e_1,\dots,\hat e_i,\dots,e_{n+1}]=0$ for
$r+1\leqslant i\leqslant n+1$.
\end{itemize}
In what follows we shall consider algebras over fields of characteristic zero, but it is
easy to see that some of our results can be extended to fields of positive characteristic.

Let $\varphi$ be a $\delta$-derivation of an algebra $A$. Then we denote by $[\varphi]$ the
matrix of the linear map $\varphi$ in the basis $\{e_1,\dots,e_{n+1}\}$. Thus,
$[\varphi(x)]=[x][\varphi]$, where $[x]$ is the row vector composed of the coordinates of the
vector $x$ in the basis $\{e_1,\dots,e_{n+1}\}$.

It is easy to see that for an algebra of type $(A_1)$ any endomorphism is a
$\delta$-derivation for any $\delta\in F$, and any endomorphism lies in $\Gamma(A)$.

\begin{lemma}\label{l1}
Let $A$ be an algebra of type $(B_1)$ and $\varphi$ a $\delta$-derivation
of $A$. Then
\[
\varphi(e_1)=\delta\sum_{k=2}^{n+1}\beta_{kk}e_1,\qquad
\varphi(e_i)=\sum_{j=1}^{n+1}\beta_{ij}e_j,\quad i\geqslant2.
\]
Conversely, every map of this form is a $\delta$-derivation. In particular, the map
$\varphi(e_2)=e_3$, $\varphi(e_i)=0$ for $i\neq2$, is a non-trivial $\delta$-derivation
of~$A$.
\end{lemma}

\begin{proof}
Let $\varphi(e_i)=\sum_{j=1}^{n+1}\beta_{ij}e_j$ with $\beta_{ij}\in F$. We note that
\[
\varphi(e_1)=\varphi[e_2,\dots,e_{n+1}]=\delta\sum_{k=2}^{n+1}[e_2,\dots,\varphi(e_k),\dots,e_{n+1}]
=\delta\sum_{k=2}^{n+1}\beta_{kk}e_1.
\]
It is easy to see that the map given by the rule
\[
\varphi(e_1)=\delta\sum_{k=2}^{n+1}\beta_{kk}e_1,\qquad
\varphi(e_i)=\sum_{j=1}^{n+1}\beta_{ij}e_j,\quad i\geqslant2,
\]
is a $\delta$-derivation. The map $\varphi(e_2)=e_3$, $\varphi(e_i)=0$ for $i\neq2$, is
of this form, and $\varphi\notin\Gamma(A)$, since $[e_2,e_2,e_4,\dots,e_{n+1}]=0$, whereas
$[e_2,\varphi(e_2),e_4,\dots,e_{n+1}]=[e_2,\dots,e_{n+1}]=e_1$.
\end{proof}

\begin{lemma}\label{l2}
Let $A$ be an algebra of type $(B_2)$ and $\varphi$ a $\delta$-derivation
of $A$. Then
\[
\varphi(e_1)=\frac{\delta}{1-\delta}\sum_{k=2}^{n}\beta_{kk}e_1,\qquad
\varphi(e_{n+1})=\beta_{n+1,n+1}e_{n+1},\qquad
\varphi(e_i)=\sum_{j=1}^{n+1}\beta_{ij}e_j,\quad 2\leqslant i\leqslant n.
\]
Conversely, every map of this form is a $\delta$-derivation. In particular, the map
$\varphi(e_2)=e_1$, $\varphi(e_i)=0$ for $i\neq2$, is a non-trivial $\delta$-derivation
of~$A$.
\end{lemma}

\begin{proof}
Let $\varphi(e_i)=\sum_{j=1}^{n+1}\beta_{ij}e_j$ with $\beta_{ij}\in F$. We observe that
\[
\varphi(e_1)=\varphi[e_1,\dots,e_n]=\delta\sum_{k=1}^{n}[e_1,\dots,\varphi(e_k),\dots,e_n]
=\delta\sum_{k=1}^{n}\beta_{kk}e_1,
\]
whence $\beta_{11}=\frac{\delta}{1-\delta}\sum_{k=2}^n\beta_{kk}$. Note that
\[
0=\varphi[e_1,\dots,\hat e_i,\dots,e_{n+1}]=\delta[e_1,\dots,\hat e_i,\dots,e_n,\beta_{n+1,i}e_i],
\]
which implies $\beta_{n+1,i}=0$ for $i\neq n+1$. It is easy to see that the map given by the
rule
\[
\varphi(e_1)=\frac{\delta}{1-\delta}\sum_{k=2}^{n}\beta_{kk}e_1,\qquad
\varphi(e_{n+1})=\beta_{n+1,n+1}e_{n+1},\qquad
\varphi(e_i)=\sum_{j=1}^{n+1}\beta_{ij}e_j,\quad 2\leqslant i\leqslant n,
\]
is a $\delta$-derivation. The map $\varphi(e_2)=e_1$,
$\varphi(e_i)=0$ for $i\neq2$, is of this form, and $\varphi\notin\Gamma(A)$, since
$[e_2,e_2,e_3,\dots,e_n]=0$, whereas
$[\varphi(e_2),e_2,e_3,\dots,e_n]=[e_1,\dots,e_n]=e_1$.
\end{proof}

\begin{lemma}\label{l3}
Let $A$ be an algebra of type $(C_1)$ and $\varphi$ a $\delta$-derivation
of $A$. Then $[\varphi]=\begin{pmatrix}A&C\\G&B\end{pmatrix}$, where $A\in M_{n-1}$,
$B\in M_2$, $C\in M_{n-1,2}$, $G=0\in M_{2,n-1}$, and
if $\delta\neq-1$, then for $w=\frac{\delta}{1-\delta}\sum_{k=1}^{n-1}\beta_{kk}$ we
have
\[
B=\begin{pmatrix}w&\beta_{n,n+1}\\ \frac{\delta}{\alpha}\beta_{n,n+1}&w\end{pmatrix},
\]
where $\beta_{n,n+1}=0$; if $\delta=-1$, then
\[
B=\begin{pmatrix}\beta_{nn}&\beta_{n,n+1}\\ -\frac{1}{\alpha}\beta_{n,n+1}&\beta_{n+1,n+1}\end{pmatrix},
\qquad \beta_{nn}+\beta_{n+1,n+1}=-\sum_{k=1}^{n-1}\beta_{kk}.
\]
Conversely, every map of this form is a $\delta$-derivation. In particular, the map
$\varphi(e_1)=e_n$, $\varphi(e_i)=0$ for $i\neq1$, is a non-trivial $\delta$-derivation
of~$A$.
\end{lemma}

\begin{proof}
Let $\varphi(e_i)=\sum_{j=1}^{n+1}\beta_{ij}e_j$ with $\beta_{ij}\in F$. For
$i\leqslant n-1$ we have
\[
0=\varphi[e_1,\dots,\hat e_i,\dots,e_{n+1}]=\pm\delta\beta_{ni}e_n\pm\delta\alpha\beta_{n+1,i}e_{n+1},
\]
whence $G=0$. We observe that
\[
\alpha\varphi(e_{n+1})=\delta\sum_{k=1}^n[e_1,\dots,\varphi(e_k),\dots,e_n]
=\alpha\delta\sum_{k=1}^n\beta_{kk}e_{n+1}+\delta\beta_{n,n+1}e_n,
\]
whence
\[
\varphi(e_{n+1})=\delta\sum_{k=1}^n\beta_{kk}e_{n+1}+\frac{\delta}{\alpha}\beta_{n,n+1}e_n.
\]
It is easy to see that
\[
\varphi(e_n)=\varphi[e_1,\dots,e_{n-1},e_{n+1}]=\delta\sum_{k=1}^{n-1}\beta_{kk}e_n
+\delta^2\sum_{k=1}^{n}\beta_{kk}e_n+\delta^2\beta_{n,n+1}e_{n+1},
\]
whence we obtain that
$(1-\delta^2)\beta_{nn}=\delta(1+\delta)\sum_{k=1}^{n-1}\beta_{kk}$. Hence, if
$\delta\neq-1$, then
{\small\[
\beta_{nn}=\frac{\delta}{1-\delta}\sum_{k=1}^{n-1}\beta_{kk},\qquad
\beta_{n+1,n+1}=\delta\sum_{k=1}^{n}\beta_{kk}
=\delta\Bigl(\sum_{k=1}^{n-1}\beta_{kk}+\frac{\delta}{1-\delta}\sum_{k=1}^{n-1}\beta_{kk}\Bigr)
=\frac{\delta}{1-\delta}\sum_{k=1}^{n-1}\beta_{kk},
\]}
and if $\delta=-1$, then $\beta_{n+1,n+1}=-\sum_{k=1}^{n}\beta_{kk}$, that is,
$\beta_{nn}+\beta_{n+1,n+1}=-\sum_{k=1}^{n-1}\beta_{kk}$. We also see that
$\beta_{n,n+1}=\delta^2\beta_{n,n+1}$, that is, $\beta_{n,n+1}\neq0$ only for $\delta=-1$.

We set $w=\frac{\delta}{1-\delta}\sum_{k=1}^{n-1}\beta_{kk}$. It is easy to see that the map
$\varphi$ defined by the formulae
\[
\varphi(e_i)=\sum_{j=1}^{n+1}\beta_{ij}e_j,\ i<n,\qquad
\varphi(e_n)=we_n+\beta_{n,n+1}e_{n+1},\qquad
\varphi(e_{n+1})=\frac{\delta}{\alpha}\beta_{n,n+1}e_n+we_{n+1},
\]
where $\beta_{n,n+1}=0$ and $\delta\neq-1$, is a
$\delta$-derivation, and for $\delta=-1$ the map defined by the formulae
 \begin{center}$\varphi(e_i)=\sum_{j=1}^{n+1}\beta_{ij}e_j$, $i<n$,
$\varphi(e_n)=\beta_{nn}e_n+\beta_{n,n+1}e_{n+1}$,
$\varphi(e_{n+1})=-\frac{1}{\alpha}\beta_{n,n+1}e_n+\beta_{n+1,n+1}e_{n+1}$,
\end{center}  where
$\beta_{nn}+\beta_{n+1,n+1}=-\sum_{k=1}^{n-1}\beta_{kk}$, is a $\delta$-derivation.
 The map
$\varphi(e_1)=e_n$, $\varphi(e_i)=0$ for $i\neq1$, is of this form, and
$\varphi\notin\Gamma(A)$, since $[e_1,\dots,e_{n-1},e_1]=0$, whereas
$[e_1,\dots,e_{n-1},\varphi(e_1)]=[e_1,\dots,e_n]=\alpha e_{n+1}\neq0$.
\end{proof}

\begin{lemma}\label{l4}
Let $A$ be an algebra of type $(C_2)$ and $\varphi$ a $\delta$-derivation
of $A$. Then $[\varphi]=\begin{pmatrix}A&C\\G&B\end{pmatrix}$, where $A\in M_{n-1}$,
$B\in M_2$, $C\in M_{n-1,2}$, $G=0\in M_{2,n-1}$, and
\[
B=\begin{pmatrix}
\frac{\delta\,\mathrm{tr}(A)}{1-\delta}-\frac{\beta\delta}{1+\delta}\gamma & \gamma\\[4pt]
\delta\gamma & \frac{\delta\,\mathrm{tr}(A)}{1-\delta}+\frac{\beta\delta}{1+\delta}\gamma
\end{pmatrix},
\]
where $\gamma=\beta_{n,n+1}$ and $\gamma$$=0$ if
$\delta\neq\frac{-\beta^2-2\pm\sqrt{\beta^4+4\beta^2}}{2}$ or $\delta=-1$; if
$\delta=-1$, then $\gamma=0$ and $\beta_{nn}=\beta_{n+1,n+1}=-\frac{\mathrm{tr}(A)}{2}$.
Conversely, every map of this form is a $\delta$-derivation. In particular, the map
$\varphi(e_1)=e_n$, $\varphi(e_i)=0$ for $i\neq1$, is a non-trivial $\delta$-derivation
of~$A$.
\end{lemma}

\begin{proof}
Let $\varphi(e_i)=\sum_{j=1}^{n+1}\beta_{ij}e_j$ with $\beta_{ij}\in F$. For
$i\leqslant n-1$ we have
\[
0=\varphi[e_1,\dots,\hat e_i,\dots,e_{n+1}]=\pm\delta\beta_{ni}(e_n+\beta e_{n+1})\pm\delta\beta_{n+1,i}e_{n+1},
\]
whence $G=0$. We observe that
\[
\varphi(e_{n+1})=\varphi[e_1,\dots,e_n]
=\delta\Bigl(\sum_{i=1}^n\beta_{ii}e_{n+1}+\beta_{n,n+1}(e_n+\beta e_{n+1})\Bigr)
=\beta_{n+1,n}e_n+\beta_{n+1,n+1}e_{n+1}.
\]
Therefore,
\begin{multline*}
\varphi(e_n)+\beta\beta_{n+1,n}e_n+\beta\beta_{n+1,n+1}e_{n+1}=\varphi(e_n+\beta e_{n+1})
=\varphi[e_1,\dots,\hat e_n,e_{n+1}]\\
=\delta\Bigl(\sum_{i=1}^{n-1}\beta_{ii}(e_n+\beta e_{n+1})+\beta_{n+1,n}e_{n+1}
+\beta_{n+1,n+1}(e_n+\beta e_{n+1})\Bigr).
\end{multline*}
Writing $\theta=\sum_{i=1}^{n-1}\beta_{ii}$, we find that
\begin{align}
\beta_{nn}&=\delta\theta+\delta\beta_{n+1,n+1}-\beta\beta_{n+1,n},\label{e4}\\
\beta_{n,n+1}&=\delta\theta\beta+\delta\beta_{n+1,n}+\delta\beta\beta_{n+1,n+1}-\beta\beta_{n+1,n+1},\label{e5}\\
\beta_{n+1,n}&=\delta\beta_{n,n+1},\label{e6}\\
\beta_{n+1,n+1}&=\delta\theta+\delta\beta_{nn}+\delta\beta\beta_{n,n+1}.\label{e7}
\end{align}
By adding and subtracting \eqref{e4} and \eqref{e7} we obtain that
\begin{align}
(1-\delta)(\beta_{nn}+\beta_{n+1,n+1})&=2\delta\theta,\label{e8}\\
(1+\delta)(\beta_{nn}-\beta_{n+1,n+1})&=-2\delta\beta\beta_{n,n+1}.\label{e9}
\end{align}
It follows easily from \eqref{e8} and \eqref{e9} that for $\delta\neq-1$ we have
\begin{align}
\beta_{nn}&=\frac{\delta\theta}{1-\delta}-\frac{\beta\delta}{1+\delta}\beta_{n,n+1},\label{e10}\\
\beta_{n+1,n+1}&=\frac{\delta\theta}{1-\delta}+\frac{\beta\delta}{1+\delta}\beta_{n,n+1},\label{e11}
\end{align}
and for $\delta=-1$ we have $\beta_{n+1,n}=\beta_{n,n+1}=0$
and, by \eqref{e8} and \eqref{e5}, $\beta_{nn}=\beta_{n+1,n+1}=-\frac{\theta}{2}$.

Substituting the expressions for $\beta_{n+1,n+1}$ and $\beta_{n+1,n}$ from
equations~\eqref{e11} and~\eqref{e6} into~\eqref{e5}, we obtain that
\[
(1-\delta^2)\beta_{n,n+1}=\frac{\beta^2\delta(\delta-1)}{1+\delta}\beta_{n,n+1},
\]
that is, $\beta_{n,n+1}\neq0$ only for
$\delta=\frac{-\beta^2-2\pm\sqrt{\beta^4+4\beta^2}}{2}$. Thus, we have shown that the map
$\varphi$ has the form described in the statement of the lemma. Clearly,
the map thus defined is a $\delta$-derivation. The map $\varphi(e_1)=e_n$, $\varphi(e_i)=0$ for $i\neq1$, is
of this form, and $\varphi\notin\Gamma(A)$, since $[e_1,\dots,e_{n-1},e_1]=0$, whereas
$[e_1,\dots,e_{n-1},\varphi(e_1)]=[e_1,\dots,e_n]=e_{n+1}$.
\end{proof}

\begin{lemma}\label{l5}
Let $A$ be an algebra of type $(D_r)$ and $\varphi$ a $\delta$-derivation
of $A$. Then $[\varphi]=\begin{pmatrix}A&C\\G&B\end{pmatrix}$, where $A\in M_r$,
$B\in M_{n+1-r}$, $C\in M_{r,n+1-r}$, $C=0$, and
furthermore,
\begin{itemize}
\item[1)] if $\delta=-1$, then $A=(a_{ij})$, where $a_{ij}=(-1)^{i-j}a_{ji}$ for $i\neq j$,
and $\mathrm{tr}(A)=-\mathrm{tr}(B)$;
\item[2)] if $\delta=\frac1{r-1}$, then $\mathrm{tr}(B)=0$ and $A=\beta E$ for $\beta\in F$;
\item[3)] if $\delta\neq\frac1{r-1},-1$, then
$A=\frac{\delta}{1+\delta-r\delta}\mathrm{tr}(B)E$;
\item[4)]
conversely, every map of this form is a $\delta$-derivation; in particular, if
$r\leqslant n$, then the map $\varphi(e_{n+1})=e_2$, $\varphi(e_i)=0$ for $i\leqslant n$, is
a non-trivial $\delta$-derivation of~$A$.
\end{itemize}
\end{lemma}

\begin{proof}
Let $\varphi(e_i)=\sum_{j=1}^{n+1}\beta_{ij}e_j$ with $\beta_{ij}\in F$. We observe that for
$i\leqslant r$ we have
\begin{multline*}
\varphi(e_i)=\varphi[e_1,\dots,\hat e_i,\dots,e_{n+1}]
=\delta\sum_{j\neq i}[e_1,\dots,\varphi(e_j),\dots,\hat e_i,\dots,e_{n+1}]\\
=\delta\sum_{j\neq i}\beta_{jj}e_i+\delta\sum_{j\leqslant r,\,j\neq i}(-1)^{i-j+1}\beta_{ji}e_j.
\end{multline*}
Hence for $j,i\leqslant r$ we have
\[
\beta_{ii}=\delta\sum_{k\neq i}\beta_{kk},\qquad
\beta_{ij}=(-1)^{i-j+1}\delta\beta_{ji}=\delta^2\beta_{ij},
\]
that is, for $i,j\leqslant r$, $i\neq j$, we have either $\beta_{ij}=(-1)^{i-j}\beta_{ji}$
and $\delta=-1$, or $\beta_{ij}=0$ and $\delta\neq-1$. We also obtain that for $i\leqslant r$
we have
\[
(1+\delta)\beta_{ii}=\delta\sum_{j=1}^{n+1}\beta_{jj},
\]
that is, either $\beta_{ii}=\beta$ and $\delta\neq-1$, or $\sum_{j=1}^{n+1}\beta_{jj}=0$ and
$\delta=-1$. Hence we see that assertion~1) holds.

If $\delta\neq-1$ and $i\leqslant r$, then
\[
\beta_{ii}=\frac{\delta}{1+\delta}\mathrm{tr}[\varphi]=\frac{r\delta}{1+\delta}\beta_{ii}
+\frac{\delta}{1+\delta}\mathrm{tr}(B),
\]
that is, either $\beta_{ii}=\frac{\delta}{1+\delta-r\delta}\mathrm{tr}(B)$ and
$\delta\neq\frac1{r-1}$, or $\mathrm{tr}(B)=0$ and $\delta=\frac1{r-1}$.

When $i>r$, standard operations yield that
\[
0=\varphi[e_1,\dots,\hat e_i,\dots,e_{n+1}]=\delta\sum_{j\leqslant r}\beta_{ji}e_j,
\]
that is, $\beta_{ji}=0$ if $i>r$ and $j\leqslant r$, and so $C=0$.

It is easy to see that a map $\varphi$ defined as in the hypothesis of the lemma is a
$\delta$-derivation.

 If $r\leqslant n$, then
the map $\varphi(e_{n+1})=e_2$, $\varphi(e_i)=0$ for $i\leqslant n$, is of this form, and
$\varphi\notin\Gamma(A)$, since $[e_{n+1},e_3,\dots,e_{n+1}]=0$, whereas
$[\varphi(e_{n+1}),e_3,\dots,e_{n+1}]=[e_2,\dots,e_{n+1}]=e_1$. Hence assertion~4) holds.
\end{proof}

\begin{theorem}\label{t6}
Every non-simple $(n+1)$-dimensional $n$-ary Filippov algebra over an algebraically closed
field of characteristic zero, except the algebra of type $(A_1)$, has a non-trivial
$\delta$-derivation for an arbitrary $\delta\neq0,1$. The algebra of type $(A_1)$ has no
non-trivial $\delta$-derivations.
\end{theorem}

\begin{proof}
In accordance with the classification of $(n+1)$-dimensional $n$-ary Filippov
algebras~\cite{b20}, all such algebras that are not simple are exhausted by the following
types: $(A_1)$, $(B_1)$, $(B_2)$, $(C_1)$, $(C_2)$, $(D_r)$ for $r\neq n+1$. For the
algebra of type $(A_1)$ every endomorphism lies in $\Gamma(A)$. Thus, the proof of the
theorem follows from the classification of algebras given above along with
Lemmas~\ref{l1}--\ref{l5}.
\end{proof}

\begin{theorem}\label{t7}
Every simple finite-dimensional Filippov algebra over an algebraically closed field of
characteristic zero has non-trivial anti-derivations and does not have non-trivial
$\delta$-derivations other than anti-derivations.
\end{theorem}

\begin{proof}
According to~\cite{b21} every simple finite-dimensional $n$-ary Filippov algebra over an
algebraically closed field of characteristic zero is isomorphic to a simple $n$-ary Filippov
algebra of type $(D_{n+1})$, as first described by Filippov~\cite{b20}. Thus, the assertion
of the theorem follows from Lemma~\ref{l5}: for $r=n+1$ we have $[\varphi]=A$; if
$\delta\neq-1$, then by assertions~2) and~3) of Lemma~\ref{l5} the map $\varphi$ is a scalar
multiple of the identity map, and hence $\varphi\in\Gamma(A)$; if $\delta=-1$, then by
assertion~1) of Lemma~\ref{l5} the map $\varphi(e_1)=e_2$, $\varphi(e_2)=-e_1$,
$\varphi(e_i)=0$ for $i\geqslant3$, is an anti-derivation, and $\varphi\notin\Gamma(A)$,
since \begin{center}$[e_1,e_1,e_4,\dots,e_{n+1}]=0$, whereas
$[e_1,\varphi(e_1),e_4,\dots,e_{n+1}]=[e_1,e_2,e_4,\dots,e_{n+1}]=e_3$.
\end{center} \end{proof}

We point out that $n$-ary Filippov algebras of type $(D_{n+1})$ generalize (to the case of
$n$-ary operations) the well-known simple Lie algebra $\mathrm{sl}_2$, and coincide with it
when $n=2$. Hopkins studied anti-derivations of simple finite-dimensional Lie algebras
in~\cite{b1}, where she constructed examples of non-trivial anti-derivations for the algebra
$\mathrm{sl}_2$. Our results agree with those in~\cite{b1} and~\cite{b3,b4} on
anti-derivations and $\delta$-derivations of the simple Lie algebra $\mathrm{sl}_2$. Note
that the derivations of simple finite-dimensional Filippov algebras over an algebraically
closed field of characteristic zero were described in~\cite{b20}.

\section{$\delta$-Derivations of the ternary Malcev  algebra $M_8$}

The class of $n$-ary Malcev  algebras was defined in~\cite{b22} as a certain natural class
of $n$-ary algebras containing the class of $n$-ary algebras of a cross product of vectors.
At present the only known example of a simple $n$-ary Malcev  algebra that is not a
Filippov algebra is the simple ternary Malcev  algebra $M_8$ arising on an 8-dimensional
composition algebra. The derivations of $M_8$ were described in~\cite{b23}, and
in~\cite{b24} its root decomposition was constructed and the structure of a
$\mathbb Z_3$-grading was introduced.

We define the \emph{$n$-ary Jacobian} to be the following function defined on an $n$-ary
algebra:
\begin{multline*}
J(x_1,\dots,x_n;y_2,\dots,y_n)=\bigl[[x_1,\dots,x_n],y_2,\dots,y_n\bigr]\\
-\sum_{i=1}^n\bigl[x_1,\dots,[x_i,y_2,\dots,y_n],\dots,x_n\bigr].
\end{multline*}
It follows from the definition that if $A$ is an $n$-ary Filippov algebra, then
\[
J(x_1,\dots,x_n;y_2,\dots,y_n)=0
\]
for all $x_1,\dots,x_n,y_2,\dots,y_n\in A$.

An \emph{$n$-ary Malcev  algebra} ($n\geqslant3$) is defined as an algebra $L$ with one
anticommutative $n$-ary operation $[x_1,\dots,x_n]$ satisfying the identity
\[
-J(zR_x,x_2,\dots,x_n;y_2,\dots,y_n)=J(z,x_2,\dots,x_n;y_2,\dots,y_n)R_x,
\]
where $R_x=R_{x_2,\dots,x_n}$ is the operator of right multiplication:
$zR_x=[z,x_2,\dots,x_n]$.

In what follows we assume that $F$ is a field of characteristic different from $2,3$ and
denote by $A$ a composition algebra over $F$ with involution $a\to\overline a$ and
identity~$1$ (see, for example,~\cite{b25}). We assume that the bilinear form
$(x,y)=\frac12(x\overline y+y\overline x)$ defined on $A$ is non-degenerate and denote by
$n(a)$ the norm of an element $a\in A$. We define a ternary multiplication operation
$[\,\cdot\,,\cdot\,,\cdot\,]$ on $A$ by the rule
\[
[x,y,z]=x\overline yz-(y,z)x+(x,z)y-(x,y)z.
\]
Then $A$ becomes a ternary Malcev  algebra~\cite{b22}, which is denoted by $M(A)$, and by
$M_8$ when $\dim(A)=8$.

Recall that a \emph{derivation of a ternary algebra} is a linear map $D$ satisfying
equation~\eqref{e2} with $n=3$. The derivations of $M_8$ were described in~\cite{b23} and it
was shown that every derivation is inner, that is,
\[
\mathrm{Der}(M_8)=\langle[R_{x,y},R_{x,z}]+R_{x,[y,x,z]}\mid x,y,z\in M_8\rangle.
\]

By a $\delta$-derivation of a ternary algebra we mean a linear map $\varphi$ satisfying
equation~\eqref{e3} with $n=3$. In the case of a ternary algebra $L$, every element of the
centroid $\Gamma(L)$ is clearly a $\frac13$-derivation. We regard a non-zero
$\delta$-derivation $\varphi$ as non-trivial if $\delta\neq0,1$ and
$\varphi\notin\Gamma(L)$.

Let $U$ be a subspace of $V$ and let $x\in V$. Let $x\big|_U$ denote the projection of the
vector $x$ onto $U$.

\begin{theorem}\label{t8}
The ternary algebra $M_8$ has no non-trivial $\delta$-derivations.
\end{theorem}

\begin{proof}
Let $1,a,b,c$ be orthonormal vectors in $A$. We choose the following basis in $A$:
\[
\{e_1=1,\ e_2=a,\ e_3=b,\ e_4=ab,\ e_5=c,\ e_6=ac,\ e_7=bc,\ e_8=abc\}.
\]

The operator $R_{x,y}$ of right multiplication is said to be \emph{regular} if in the
Fitting decomposition $M=M_0\oplus M_1$ with respect to $R_{x,y}$ the dimension of $M_0$ is
minimal~\cite{b24}. By Theorem~1 in~\cite{b24} we have the following root decomposition of
$M_8$: $M=M_0\oplus M_\alpha\oplus M_{-\alpha}$, where $\alpha\in F$ is such that
$vR_{x,y}=\pm\alpha v$ for any $v\in M_{\pm\alpha}$. It is known from Lemma~3 in~\cite{b24}
that there exists a non-trivial grading on $M_8$. If we denote $M_{\pm\alpha}$ by $M_{\pm1}$,
then
\[
[M_i,M_j,M_k]\subseteq M_{i+j+k\,(\mathrm{mod}\,3)}.
\]

Let $v\in M_\alpha$. Then
\begin{equation}\label{e12}
\alpha\varphi(v)=\varphi(vR_{x,y})=\delta\bigl([\varphi(v),x,y]+[v,\varphi(x),y]+[v,x,\varphi(y)]\bigr).
\end{equation}

We assume that
\[
\varphi(v)\big|_{M_\alpha}=w_\alpha,\qquad
\varphi(x)\big|_{M_0}=\alpha_xx+\beta_xy,\qquad
\varphi(y)\big|_{M_0}=\alpha_yx+\beta_yy.
\]
In view of these relations, equation~\eqref{e12} yields that
$w_\alpha=\frac{\alpha_x+\beta_y}{1-\delta}\delta v$. We observe that by Lemma~1
in~\cite{b24} the operators $R_{e_i,e_i+e_j}$ and $R_{e_j,e_i+e_j}$ for $i\neq j$ are regular.
Consequently, we can conclude that
\[
\frac{\alpha_{e_i}+\beta_{e_i+e_j}}{1-\delta}\delta v=\frac{\alpha_{e_j}+\beta_{e_i+e_j}}{1-\delta}\delta v,
\]
whence $\alpha_{e_i}=\alpha_{e_j}$. Since $\frac{2\alpha_{e_i}}{1-\delta}\delta=\alpha_{e_i}$,
we obtain that either $\delta=\frac13$ or $\alpha_{e_i}=0$.

For every $i\in\{2,\dots,8\}$ it is possible to choose $j,k,l,m,s,t$ depending on $i$ such
that
\begin{align}
e_i&=e_je_k=e_le_m=e_se_t,\label{e13}\\
e_j&=e_se_m=e_ke_i=e_te_l,\label{e14}\\
e_k&=e_ie_j=e_me_t=e_se_l,\label{e15}\\
e_l&=e_me_i=e_ke_s=e_je_t,\label{e16}\\
e_m&=e_ie_l=e_te_k=e_je_s,\label{e17}\\
e_s&=e_le_k=e_te_i=e_me_j,\label{e18}\\
e_t&=e_ie_s=e_ke_m=e_le_j.\label{e19}
\end{align}

We can assume that $\varphi(e_q)=\sum_{p=1}^8a_{qp}e_p$. We already know that
$a_{qq}=a_{pp}$ and that $a_{pp}=0$ when $\delta\neq\frac13$. Since
\[
\varphi(e_s)=\varphi([e_i,e_j,e_l])
=\delta\bigl([\varphi(e_i),e_j,e_l]+[e_i,\varphi(e_j),e_l]+[e_i,e_j,\varphi(e_l)]\bigr),
\]
the corresponding multiplication operations yield that
\begin{multline*}
\sum_{p=1}^8a_{sp}e_p=\delta(-a_{it}e_1+a_{ik}e_m-a_{im}e_k+a_{i1}e_t-a_{is}e_i-a_{jm}e_1
+a_{j1}e_m\\
+a_{jt}e_k-a_{jk}e_t-a_{js}e_j+a_{lk}e_1+a_{lt}e_m-a_{l1}e_k-a_{lm}e_t-a_{ls}e_l).
\end{multline*}
(For example,
$[e_i,e_l,e_k]=(e_i\overline{e_l})e_k=(e_le_i)e_k=-e_me_k=e_ke_m=e_t$.) Consequently,
$a_{sp}=-\delta a_{ps}$, where $p\in\{i,j,l\}$. Since the index $i$ is arbitrary, the
relations \eqref{e13}--\eqref{e19} enable us to conclude that $a_{pq}=-\delta a_{qp}$ for all
$p,q\in\{1,\dots,8\}$. Using this relation we can conclude that
$a_{pq}=-\delta a_{qp}=\delta^2a_{pq}$, whence we have $\delta=-1$ and $a_{pq}=a_{qp}$ or
$\varphi$ is trivial. We can now show that $M_8$ has no non-trivial anti-derivations. In
what follows we follow the scheme of argument used above and take into account the fact
that $\delta=-1$.

Suppose that $-\varphi(e_t)=\varphi([e_i,e_k,e_l])$. Then
\begin{multline*}
-\sum_{p=1}^8a_{tp}e_p=-(a_{i1}e_s-a_{ij}e_m+a_{im}e_j-a_{is}e_1+a_{it}e_i+a_{k1}e_m+a_{kj}e_s\\
-a_{km}e_1-a_{ks}e_j+a_{kt}e_k+a_{l1}e_j-a_{lj}e_1-a_{lm}e_s+a_{ls}e_m+a_{lt}e_l),
\end{multline*}
whence we obtain that
\begin{equation}\label{e20}
a_{ts}=a_{i1}+a_{jk}-a_{ml}.
\end{equation}

Suppose that $\varphi(e_m)=\varphi([e_i,e_j,e_t])$. Then
\begin{multline*}
\sum_{p=1}^8a_{mp}e_p=-(a_{i1}e_l-a_{ik}e_s-a_{il}e_1+a_{im}e_i-a_{is}e_k+a_{j1}e_s-a_{jk}e_l\\
+a_{jl}e_k+a_{jm}e_j-a_{js}e_1+a_{t1}e_k-a_{tk}e_1-a_{tl}e_s+a_{tm}e_t+a_{ts}e_l),
\end{multline*}
whence we obtain that
\begin{equation}\label{e21}
a_{ml}=a_{i1}-a_{jk}+a_{ts}.
\end{equation}

Suppose that $-\varphi(e_t)=\varphi([e_i,e_j,e_m])$. Then
\begin{multline*}
-\sum_{p=1}^8a_{tp}e_p=-(a_{i1}e_s-a_{ik}e_l+a_{il}e_k-a_{is}e_1+a_{it}e_i-a_{j1}e_l-a_{jk}e_s\\
+a_{jl}e_1+a_{js}e_k+a_{jt}e_j-a_{m1}e_k+a_{mk}e_1+a_{ml}e_s-a_{ms}e_l+a_{mt}e_m),
\end{multline*}
whence we obtain that
\begin{equation}\label{e22}
a_{ts}=a_{i1}-a_{jk}+a_{ml}.
\end{equation}

Suppose that $-\varphi(e_s)=\varphi([e_i,e_k,e_m])$. Then
\begin{multline*}
-\sum_{p=1}^8a_{sp}e_p=-(-a_{i1}e_t+a_{ij}e_l+a_{is}e_i+a_{it}e_1-a_{il}e_j-a_{k1}e_l-a_{kj}e_t\\
+a_{ks}e_k+a_{kt}e_j+a_{kl}e_1+a_{m1}e_j-a_{mj}e_1+a_{ms}e_m+a_{mt}e_l-a_{ml}e_t),
\end{multline*}
whence we obtain that
\begin{equation}\label{e23}
a_{st}=-a_{i1}-a_{jk}-a_{ml}.
\end{equation}

We observe that equations \eqref{e20} and \eqref{e21} imply that
\begin{equation}\label{e24}
a_{kj}=a_{st},\qquad a_{ml}=a_{i1}.
\end{equation}
The formulae \eqref{e20} and \eqref{e22} imply that
\begin{equation}\label{e25}
a_{i1}=a_{st},\qquad a_{ml}=a_{kj},
\end{equation}
and the relations \eqref{e23} and \eqref{e24} imply that
\begin{equation}\label{e26}
a_{st}=-a_{ml}.
\end{equation}

By analysing equations \eqref{e24}--\eqref{e26} we obtain that
\[
a_{i1}=a_{st}=a_{ml}=a_{kj}=0.
\]
Since the index $i$ is arbitrary, we obtain from the relations \eqref{e13}--\eqref{e19} that
the map $\varphi$ is trivial.
\end{proof}

Theorem~\ref{t8} implies the existence of a simple $n$-ary Malcev  algebra over an
algebraically closed field of characteristic zero that is not an $n$-ary Filippov algebra
and has no non-trivial $\delta$-derivations, in contrast to the case of $n$-ary Filippov
algebras over an algebraically closed field of characteristic zero. We point out that these
results agree with those of Filippov~\cite{b5} on the absence of non-trivial
$\delta$-derivations on prime binary non-Lie Malcev  algebras.

We note that a partial generalization of the results of the present paper was announced
in~\cite{b26}, and a detailed survey of the main results in the theory of
$\delta$-differential algebras and superalgebras was given in~\cite{b27}.

\section*{Appendix A. Differences from the published version}
\addcontentsline{toc}{section}{Appendix A. Differences from the published version}

\emph{Lemmas~\ref{l1}--\ref{l5}: statements changed.} 

\noindent\emph{Theorem~\ref{t6}: statement changed.}  

\noindent\emph{Theorem~\ref{t7}: proof completed.}

\end{document}